\documentclass{amsart}
\usepackage[a4paper, total={7in, 8in}]{geometry}

\usepackage{xcolor,colortbl}

\usepackage{algpseudocode}

\definecolor{defcolor}{HTML}{ff725e}

\usepackage{thm-restate}

\usepackage[maxbibnames=99, style = alphabetic]{biblatex}

\usepackage[]{amsfonts}
\usepackage{color}
\usepackage{comment}
\usepackage[]{amsmath}
\usepackage{amssymb}
\usepackage[]{amsthm}
\usepackage{graphicx}
\usepackage{enumerate}
\usepackage{hyperref}
\usepackage{thmtools}
\usepackage{tikz}

\title{A Proof of a conjecture of Watanabe--Yoshida via Ehrhart theory }

\author[Yakob Kahane]{Yakob Kahane}
\address[Y.~Kahane]{LACIM, Université du Québec à Montréal, Canada}
\email{kahane.yakob@courrier.uqam.ca}

\newtheorem{theorem}{Theorem}[section]
\newtheorem{lemma}[theorem]{Lemma}
\newtheorem{corollary}[theorem]{Corollary}
\newtheorem{conjecture}[theorem]{Conjecture}
\newtheorem{proposition}[theorem]{Proposition}

\newcommand{\new}[1]{\textit{\textbf{\color{blue}{#1}}}}

\theoremstyle{definition}
\newtheorem{definition}[theorem]{Definition}
\newtheorem{remark}[theorem]{Remark}
\newtheorem{example}[theorem]{Example}

\numberwithin{equation}{section}

\begin{document}

\begin{abstract}
In 2005, Watanabe and Yoshida formulated a conjecture for a lower bound of the Hilbert-Kunz multiplicity
of local rings that was recently settled by Meng using analytic methods. More recently, Pak-Shapiro-Smirnov-Yoshida
used Ehrhart theory to compute explicitly the multiplicity and reduced the conjecture to showing an inequality of
the values of the Ehrhart polynomial of a zigzag poset shifted to $t - 1/2$. We completely realize their approach to
give another proof of this Watanabe--Yoshida conjecture. The main ingredient of the proof relies on a new explicit
combinatorial formula for the coefficients of this shifted Ehrhart polynomial. In terms of the generating function of
the shifted polynomial, this formula manifests itself as a Hadamard product of the exponential generating function of
Euler numbers and an explicit algebraic function.
\end{abstract}

\maketitle 

\section{Introduction}
Let $P$ be a poset of size $n$. We are interested in the coefficients of the order polynomial  
\[
\Omega(P;t)=\sum_{k=1}^{n} c_k(P)\, t^k .
\]
This polynomial encodes substantial information about $P$. In particular, the leading coefficient of $n!\cdot \Omega(P;t)$ is the number of linear extensions of $P$, denoted $e(P)$ (see \cite[Ch.~3]{EC2}), a central and notoriously difficult to compute invariant \cite{brigthwellwrinckler}.  

From a geometric viewpoint, following Stanley \cite{stanley1986two}, we define the \emph{order polytope}
\[
\mathcal{O}(P):=\{ f\in [0,1]^P : f(u)\le f(v)\ \text{whenever }u\preceq v \text{ in } P\}.
\]
The Ehrhart polynomial of $\mathcal{O}(P)$ is equal to $\Omega(P;t+1)$. It is well known that the coefficients $n!\, c_k(P)$ are integers, though they are not always nonnegative. Determining when these coefficients are positive remains an open problem \cite{liu}. A recent result of Ferroni, Morales, and Panova \cite{ferroni2025skew} shows that to prove nonnegativity of all coefficients, it suffices to establish positivity of the linear coefficients of all convex subposets of~$P$ which are partitions of the elements of the poset (see Definition~\ref{def:2}).

\vspace{5mm}

Our starting point is a positive formula for each coefficient of the order polynomial, expressed solely in terms of the linear coefficients of convex subposets of $P$. This formula already appears in \cite[Prop.~2.9]{shareshian2003newapproachorderpolynomials}.

\begin{theorem}[\cite{shareshian2003newapproachorderpolynomials}]
Let $P$ be a poset. Then
\[
c_k(P)=\frac{1}{k!}\!
\sum_{\substack{\mu\in \mathsf{Decompo}(P)\\ \mu=\{P_1,\dots,P_k\}}}
e(\mu)\,\prod_{i=1}^k c_1(P_i),
\]
where $\mathsf{Decompo}(P)$ denotes the set of all decompositions of $P$ into $k$ convex subposets, $e(\mu)$ is the number of ways to linearly order the blocks of the decomposition, and $c_1(P_i)$ is the linear coefficient of the order polynomial of $P_i$.
\end{theorem}

This formula is particularly useful for families of posets that are closed under taking convex subposets. Zigzag posets form such a family. These posets appear in various areas of combinatorics \cite{brunner2024triangulationflowpolytopezigzag,morier-genoud-ovsienko}, as well as in algebra, geometry, and probability \cite{StanleySurveyAltPerms,chan2024correlation,coons2020hpolynomialorderpolytopezigzag,zbMATH07220122}.  

Recently, Pak, Shapiro, Smirnov, and Yoshida \cite{pak2025hilbertkunzmultiplicityquadricsehrhart} proved that the \emph{Hilbert--Kunz multiplicity} of the simple $(A_1)$-singularities of dimension $n$ and characteristic $p$, denoted $e_{HK}(\mathcal{A}_{p,n})$, satisfies an inequality that can be extracted from \cite[Cor.~3.10]{pak2025hilbertkunzmultiplicityquadricsehrhart}:

\begin{theorem}[\cite{pak2025hilbertkunzmultiplicityquadricsehrhart}]
For all $p,n$,
\begin{equation}\label{eq: HK to Order poly}
e_{HK}(\mathcal{A}_{p,n}) \ge 1 + \frac{\Omega(Z_n; p/2 - 1/2)}{(p/2)^n},
\end{equation}
where $Z_n$ is the zigzag poset of size $n$.
\end{theorem}

This raises the question of whether Ehrhart theory can be used to prove the third part of the Watanabe--Yoshida conjecture (2005). Recently, Meng announced a proof using different analytic methods \cite{meng2025analysishilbertkunztheory}.

\begin{conjecture}[\cite{watanabe2005hilbert}] \label{conj:WY}
Let $p$ be a prime and  
\[
\mathcal{A}_{p,n}:=\frac{\mathbb{F}_p[[x_0,\dots,x_n]]}{(x_1^2+\cdots+x_n^2)}.
\]
Then
\begin{equation}\label{eq: HK bigger than Ed}
e_{HK}(\mathcal{A}_{p,n})\ge 1+\frac{E_n}{n!},
\end{equation}
where $E_n$ is the Euler number, equivalently the number of linear extensions of the zigzag poset of size $n$.
\end{conjecture}

In this paper, we provide an explicit formula for each coefficient of $\Omega(Z_n;t-1/2)$, which we call the \emph{shifted Ehrhart polynomial} of the zigzag poset of size $n$.

\begin{theorem}[(Theorem~\ref{thm:3})]
Let $n\ge 2$ and $1\le k\le n$. Then
\[
[t^k]\:\Omega(Z_n;t-1/2)=\frac{n!}{k!}E_k\!\sum_{\mu\in \mathrm{Decompo}_k(Z_n)} w(\mu),
\]
where $w(\mu)$ is a weight function described in detail in Section~\ref{sec:55}.
\end{theorem}

Several recursive formulas for the ordinary generating function of $\Omega(Z_n;t)$ are known \cite[Thm.~3.15]{petersen-zhuang}, \cite[Ex.~3.66]{EC2}. We give instead a formula for the \emph{exponential} generating function of the shifted order polynomials of zigzag posets, expressed as a Hadamard product $\odot$ (in $t$) of the generating function of Euler numbers with an explicit algebraic function. To our knowledge, this is the first non-recursive expression for this generating function.

\begin{theorem}\label{thm: Hadamard order poly zigzag}
For all $t$ and $n\ge 1$,
\[
\sum_{n\ge 0} \Omega(Z_n;t-1/2)\frac{z^n}{n!}
= \bigl(\tan t+\sec t\bigr)\ \odot\ G(z;t),
\]
where  
\[
G(z;x)=\sqrt{1-z^2/4}\,\frac{1+2x\,\arcsin(z/2)}{1-4x^2\,\arcsin(z/2)}.
\]
\end{theorem}

Using this explicit formula, we obtain the following inequality.

\begin{corollary}[Corollary~\ref{corro:finalinequality}]
For $n\ge 8$ and $t\ge 3/2$,
\[
\Omega(Z_n;t-1/2) \ge t^n\frac{E_n}{n!}.
\]
\end{corollary}

By computing the remaining cases of the Hilbert--Kunz multiplicity using the formula of Pak, Shapiro, Smirnov, and Yoshida, we verify that the third part of the Watanabe--Yoshida conjecture holds for all $n\ge 2$ and $p\ge 3$. This yields an Ehrhart-theoretic proof of the conjecture, complementing the recent analytic proof announced by Meng \cite{meng2025analysishilbertkunztheory}.

\begin{theorem}[Theorem~\ref{thm:proof 3rd part WS conj}]
Conjecture~\ref{conj:WY} holds.
\end{theorem}


\section*{Outline}

In Section~3, following the work of Ferroni, Morales, and Panova \cite{ferroni2025skew}, provides a combinatorial interpretation of the coefficients of the Ehrhart polynomial of a poset polytope. This is summarized in Theorem~\ref{thm:1}. In Section~4 gives an explicit formula for the coefficients $f_{n,k}$, the coefficients of $\Omega(Z_n;t-\tfrac{1}{2})$. This is summarized in Theorem~\ref{thm:3}. In Section~5, we show that these coefficients can be computed using the Taylor expansion of an algebraic function $G(z,x)$. In Section~6, we apply the results from the previous section to obtain a lower bound for $\Omega(Z_n;t-\tfrac{1}{2})$, thereby proving the third part of the Watanabe--Yoshida conjecture.In Section~7, we discuss some final remarks that arise from our results.

\section{Notation and background}

\subsection{Notation}

In this paper, we denote by $[f]_{t^n}$ the coefficient of $t^n$ in the Taylor expansion of a function $f$, i.e :

\begin{definition}
    \[
        [f]_{t^n} := [t^n]\; f.
    \]
\end{definition}

Given a function $f$, we also define its even and odd parts:

\begin{definition}
    \[
        f_{\text{even}}(x) := \frac{f(x) + f(-x)}{2}, \qquad 
        f_{\text{odd}}(x) := \frac{f(x) - f(-x)}{2}.
    \]
\end{definition}

\subsection{Linear extension, Ideal}
We remind standards objects related to posets.
\begin{definition}
    Let $P$ be a poset of size $n$, a \new{linear extension} is a surjective increasing function from $P$ to $\{1,\dots, n\}$. The set of linear extension is denoted \new{$\mathcal{L}(P)$}. And we denote \new{$e(P)$} $:= |\mathcal{L}(P)|$. 
\end{definition}

\begin{definition}
    Let $P$ be a poset, an \new{ideal} $I$ of $P$ is a subposet of  $P$ with the following conditions.
    \begin{itemize}
        \item $I$ is not empty
        \item for every $x \in I$, and $y \in P$, $y \leq x$ implies $y \in I$. 
    \end{itemize}
\end{definition}

\subsection{Order Polynomial and Ehrhart polynomial}

\begin{definition}
    Let $P = (X,\preceq)$ be a partially ordered set (poset) with $|X| = n$.  
    The \new{order polynomial} of $P$ is the unique polynomial \new{$\Omega(P;t)$} $\in \mathbb{Q}[t]$ such that
    \begin{equation}
        \Omega(P;t) := \#\{ f : X \to [1,\dots,t] \;\text{ such that }\; f(x) \le f(y) \text{ whenever } x \preceq y \}.
    \end{equation}
\end{definition}

It is a classical result (see for example \cite{stanley2011enumerative}) that this counting function is a polynomial and that  
$n!\cdot \Omega(P;t)$ has integer coefficients. Thus, we may write:
\[
    \sum_{k=0}^n a_k(P)t^k := n!\,\Omega(P;t),
    \qquad
    \sum_{k=0}^n c_k(P)t^k := \Omega(P;t).
\]

\begin{definition}
    The \new{order polytope} of $P$, introduced by Stanley \cite{stanley1986two}, is defined as
    \begin{equation}
        \mathcal{O}(P) := \{ f \in [0,1]^P \;\text{ such that }\; f(u) \le f(v) \text{ whenever } u \preceq v \}.
    \end{equation}
\end{definition}

The polytope $\mathcal{O}(P)$ is integral. Stanley showed in \cite[Thm.~4.1]{stanley1986two} that the \new{Ehrhart polynomial}  of $\mathcal{O}(P)$, which counts the lattice points in the dilation $t\cdot \mathcal{O}(P)$, is given by the order polynomial.

\begin{proposition}[{\cite[]{stanley1986two}}]
    For any poset $P$, the Ehrhart polynomial of $\mathcal{O}(P)$ satisfies
    \[
        \#(t\cdot \mathcal{O}(P) \cap \mathbb{Z}^n) = \Omega_P(t+1).
    \]
\end{proposition}

\vspace{5mm}

In \cite{kreweras1965classe}, Kreweras gave the following determinant identity for $\Omega(P_{\lambda/\mu}; t)$:

\begin{theorem}[Kreweras]
    \begin{equation}
        \Omega(P_{\lambda/\mu}; t)
        = \det\!\left[\, \binom{\lambda_i - \mu_j + t - 1}{\lambda_i - \mu_j - i + j} \right]_{i,j=1}^{\ell},
    \end{equation}
where $\ell$ is the length of $\lambda$.
\end{theorem}

\vspace{5mm}

Recently, Ferroni, Morales, and Panova \cite[Prop.~3.3]{ferroni2025skew} showed that the coefficients $c_k$ satisfy the following recurrence:

\begin{lemma}[\cite{ferroni2025skew}]\label{lemma:1}
    For any poset $P$ and every $k \ge 2$,
    \[
        c_k(P)
        = \frac{1}{2^k - 2}
            \left(
                \sum_{\substack{I \subsetneq P \\ I \neq \emptyset}}
                \;\sum_{i=1}^{k-1} c_i(I)\,c_{k-i}(P\setminus I)
            \right),
    \]
    where the first sum is over ideals $I$ of the poset $P$.
\end{lemma}

They used this recurrence to show that the non-negativity of the coefficients of the order polynomial of a poset $P$ follows from the positivity of the linear coefficient of the order polynomial of every convex subposet of $P$.

\begin{definition}
    The \new{zeta polynomial} of a poset $P$ is a polynomial $Z(P;t) \in \mathbb{Q}[t]$ that enumerates multichains of element in $P$ of size $t$. 
\end{definition}
It is a standard fact that the order polynomial of  a poset $P$ equals the zeta polynomial of its lattice of order ideals $J(P)$. 

\vspace{5mm}

\subsection{Zigzag posets and Fibonacci Polytope}

The main poset analyzed in this article is the zigzag poset which is also the poset associated to the Fibonacci polytope. It is important in our arguments to keep the following convention on the labeling of this poset. 

\begin{definition}\label{def:5}
    The zigzag poset \new{\( Z_n \)} on the ground set \( \{z_1, \dots, z_n\} \) is the poset with cover relations
    \[
        z_1 > z_2 < z_3 > z_4 < \dots
    \]
    that is, it satisfies
    \[
        z_{2i-1} > z_{2i}
        \qquad\text{and}\qquad
        z_{2i} < z_{2i+1}
    \]
    for all \( i \) with \( 1 \le i \le \left\lfloor \frac{n-1}{2} \right\rfloor \).
\end{definition}

Using Kreweras’ formula, Ferroni, Morales, and Panova \cite{ferroni2025skew} computed the linear coefficient $c_1$ for zigzag posets:

\begin{proposition}[Ferroni--Morales--Panova {\cite[\S 5.1]{ferroni2025skew}}]
    \begin{equation}
        n!\,c_1(Z_n) =
        \begin{cases}
            \left( \frac{n-1}{2}! \right)^2, & \text{if $n$ is odd}, \\[2mm]
            \left(\frac{n}{2}\right)!\left(\frac{n-2}{2}\right)!, & \text{if $n$ is even}.
        \end{cases}
    \end{equation}
\end{proposition}

These numbers have been studied (see \cite[\href{https://oeis.org/A10551}{A010551}]{oeis}).  
The EGF is known; here we are interested only in its odd part.

\begin{proposition}
    \begin{equation}
        \sum_{k \ge 0} c_1(Z_{2k+1}) z^{2k+1}
        = \frac{2\arcsin(z/2)}{\sqrt{1 - z^2/4}}.
    \end{equation}
\end{proposition}

\begin{definition}[\cite{pak2025hilbertkunzmultiplicityquadricsehrhart}]
    The \new{$n$-dimensional Fibonacci polytope} is the subset of $[0,1]^n$ given by the inequalities
    \[
        x_i + x_{i+1} \le 1, \qquad i = 1,\dots,n-1,
    \]
    and the \new{extended $n$-dimensional Fibonacci polytope} is the subset of $[-1,1]^n$ given by
    \[
        |x_i| + |x_{i+1}| \le 1, \qquad i = 1,\dots,n-1.
    \]

    The $n$-dimensional Fibonacci polytope is the order polytope of $Z_n$.

    Let \new{$\mathsf{Fib}_n(t)$} and \new{$\mathsf{EFib}_n(t)$} denote respectively the Ehrhart polynomials of the Fibonacci polytope and the extended Fibonacci polytope of order $n$.
\end{definition}

\subsection{Euler numbers}

The Euler numbers \new{$E_n$} can be defined as the number of linear extensions of the zigzag poset:
\begin{equation}
    E_n := e(Z_n) = \#\mathcal{L}(Z_n).
\end{equation}

These numbers have the following well-known exponential generating function:

\begin{proposition}
    \begin{equation}
        \sum_{n \ge 0} \frac{E_n}{n!} z^n = \tan(z) + \sec(z).
    \end{equation}
\end{proposition}

Moreover, $E_n$ has the following explicit well-known formula (see e.g. \cite{stanley2009survey}), which will occurs later:

\begin{proposition}
    \begin{equation}
        \frac{E_n}{n!}
        = 2\left( \frac{2}{\pi} \right)^{n+1}
          \sum_{k \ge 0} (-1)^{k(n+1)} \frac{1}{(2k+1)^{n+1}}.
    \end{equation}
\end{proposition}

\subsection{Watanabe--Yoshida conjecture}
In this chapter, we do a rapid exposition of the Watanabe--Yoshida conjecture and its link with zigzgag posets. 

\begin{definition}[Monsky \cite{monsky1983hilbert}]\label{def:4}
    Let $(R,\mathfrak{m})$ be a local ring of characteristic $p>0$.  
    The \textit{Hilbert–Kunz multiplicity} \new{$e_{HK}(R)$} is the limit
    \[
        e_{HK}(R)
        := \lim_{e \to \infty}
            \frac{\lambda(R/\mathfrak{m}^{[p^e]})}{p^{e\dim R}},
    \]
    where $\mathfrak{m}^{[p^e]}$ is the ideal generated by all $p^e$-th powers of elements of $\mathfrak{m}$.
\end{definition}

Watanabe and Yoshida made the following conjecture:

\begin{conjecture}[\cite{watanabe2005hilbert}] \label{conj:1}
    Let $p$ be a prime number and define
    \[
        \textcolor{blue}{\mathcal{A}_{p,n}} :=
        \frac{\mathbb{F}_p[[x_0,\dots,x_n]]}{(x_1^2 + \cdots + x_n^2)}.
    \]
    Then
    \[
        e_{HK}(\mathcal{A}_{p,n}) \ge 1 + \frac{E_n}{n!}.
    \]
\end{conjecture}

\vspace{5mm}

Recently, Pak, Shapiro, Smirnov, and Yoshida \cite{pak2025hilbertkunzmultiplicityquadricsehrhart} obtained the following explicit formula for the Hilbert–Kunz multiplicity of $\mathcal{A}_{p,n}$.

\begin{theorem}[\cite{pak2025hilbertkunzmultiplicityquadricsehrhart}] \label{thm:2}
    \begin{equation}\label{eq:pakformula}
        e_{HK}(\mathcal{A}_{p,n})
        = 1
        + \frac{2^n \cdot \mathsf{Fib}_n\!\left(\frac{p-3}{2}\right)}
               {p^n - \mathsf{EFib}_{n-2}\!\left(\frac{p-1}{2}\right)}.
    \end{equation}
\end{theorem}

Using that $\mathsf{Fib}_n(t)$ is the shifted order polynomial of $Z_n$, and that  
$\mathsf{EFib}_{n-2}(t) \ge 0$ is a nonnegative integer whenever $t$ is a positive integer, one obtains:

\begin{equation}\label{eq:ineqwata}
    e_{HK}(\mathcal{A}_{p,n})
    \ge 1 + \frac{2^n\,\Omega(Z_n; p/2 - 1/2)}{p^n}.
\end{equation}

\section{Decomposition formula} 

We start by defining the notion of decomposition into convex subposets. It will be the cornerstone for all the results in the next sections. The aim of this section is to give a formula of each coefficient $c_k(P)$ in terms of the linear coefficients of the convex subposets of $P$. 

\begin{definition}[convex subposets]~\\
	Let $P = (X_P, \preceq_P)$ be a poset.  
A \new{convex subposet} $R$ of $P$ is a poset $R = (X_R, \preceq_R)$ satisfying the following properties:
\begin{itemize}
    \item $R$ is connected, 
    \item $X_R \subset X_P$,
    \item (Subposet) For all $a, b \in X_R$, $a \preceq_R b$ if and only if $a \preceq_P b$,
    \item (Convexity) If $a \preceq_P b \preceq_P c$ and $a, c \in X_R$, then $b \in X_R$.
\end{itemize}
\end{definition}
\vspace{5mm}

\begin{definition}[decomposition into convex subposets] \label{def:2}
	Let $P$ be a poset, a \new{decomposition $\mu$ into convex subposets} is a partition of the ground set of $P$ into $k$ convex subposets $P_1,..., P_k$.
	Let $P_i, P_j$ two subposets of $P$ from the decomposition $\mu$, we say that $P_i \tilde{\leq_{\mu}} P_j$ if and only if there exist $a_i \in P_i$ and $a_j \in P_j$ such that $a_i \leq_P a_j$. 
\end{definition}

\begin{figure}[h]
    \centering
    \includegraphics[width=0.7\linewidth]{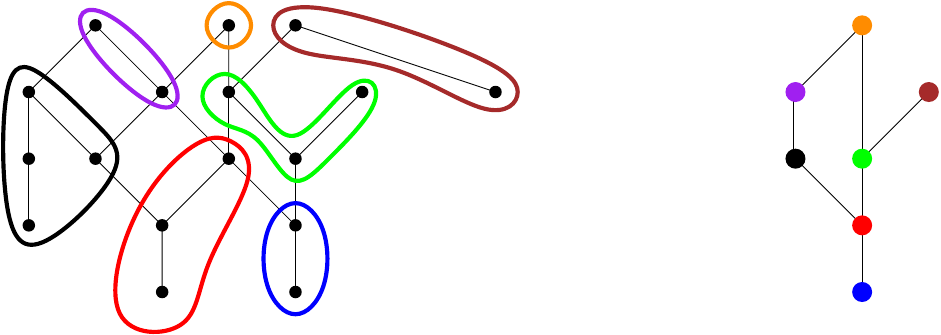}
    \caption{A decomposition $\mu$ and the induced poset.}
    \label{fig:placeholdersfds}
\end{figure}
\vspace{5mm}
Let us also define $\binom{P}{\mu} := \binom{|P|}{|P_1|, ... |P_k| }$

From the definition of subposets, we see that $P_i \tilde{\leq_{\mu}} P_j$ if and only if $i=j$. However the relation $\tilde{\leq_{\mu}}$ is not a partial order relation yet. In order to create such a relation we need to consider a smaller class of decomposition into convex subposets. 

\begin{definition}[Standard Decomposition]\label{def:3}
	A decomposition $\mu = \{P_1,...., P_k\}$ into convex subposets is said to be \new{standard} if there is no $2 \leq r \leq k$ and $r$ distinct $i_1, i_2,..., i_r$ such  that  $P_{i_1} \widetilde{\leq}_{\mu} P_{i_2} .... \widetilde{\leq}_{\mu} P_{i_r} \widetilde{\leq}_{\mu} P_{i_1}$. 
\end{definition}

\begin{proposition}\label{lemma:2}~\\
	Let $\leq_{\mu}$ be the transitive closure of $\widetilde{\leq}_{\mu}$ then : 
	\[\leq_{\mu} \text{ is a a partial order relation } \text{ if and only if } \mu \text{ is a Standard decomposition}\]
\end{proposition}
We define \new{$\text{SDecompo}_{k}(P)$} the set of standard decompositions of $P$ into $k$ convex subposets.
In this case, we define \new{$e(\mu)$} to be the number of linear extension on ${(\mu, \leq_{\mu})}$.  
Notice that a standard decomposition and one of its labelling is equivalent to a strictly increasing chains of ideals. More precisely, we have the following. 
\vspace{5mm}

\begin{proposition}\label{lemma:3}Let $P$ a poset. 
	\begin{enumerate} 
		\item  Let $\emptyset = I_0 \subset I_1 \subset .... \subset I_k \subset I_{k+1} = P$ an increasing sequence of ideals of $P$. Then $\{(I_1 \setminus I_0), (I_2 \setminus I_1),...., (I_{k+1} \setminus I_k)\}$ is a standard decomposition of $P$. 
		\item Let $\mu$ a standard decomposition, then, $e(\mu)$ is the number of increasing sequences of Ideals of $P$ such that : \\
		${\mu = \{(I_1 \setminus I_0), (I_2 \setminus I_1),...., (I_{k+1} \setminus I_k)\}}$. 
	
	\end{enumerate}
\end{proposition}
\vspace{5mm}
We can now state the main tool of this article.

\begin{restatable}{theorem}{bigrec}
\label{thm:1}
	Let $P$ a poset, then for all $k \in \mathbb{N}$ : 
    let, $\Omega(P;t)= \sum_{k=0}^d c_k(P) t^k $ the order polynomial of $P$, then 
	\[c_k(P) = \dfrac{1}{k!}\sum_{\substack{\mu \in \text{SDecompo}_{k}(P) \\  \mu = \{P_1,P_2,...,P_k\}}} e(\mu) \Pi_{i = 1}^k c_1(P_i) .\]
\end{restatable}
\begin{proof}
Let $Z([u,v];t)$ be the zeta polynomial
 of the interval $[u,v] \subseteq  P$, and let $x_1,..., x_k$ $k$ variables, we have : 
 \[Z([u,v];x_1 + \cdots + x_k) = \sum_{u = v_0 \leq v_1 \leq... \leq v_{n-1} \leq v_k = v} \prod_{i=1}^{k} Z([v_{i-1}, v_{i}],x_i).\]
 It is a standard fact that the order polynomial of a poset P equals the zeta polynomial of its lattice of order ideals $J(P)$ (see e.g \cite{stanley2011enumerative}[ch.3 p.292]).  
  \[\Omega(P;kt) = \sum_{\substack{\emptyset = v_0 \leq ... \leq v_k =P \\ v_i \in J(P)}} \Omega(v_{i+1} \setminus v_i; t).\]

 Also : 
 \[\dfrac{1}{k!}\sum_{\substack{\mu \in \mathsf{SDecompo}_{k}(P) \\  \mu = \{P_1,P_2,...,P_k\}}} e(\mu) \Pi_{i = 1}^k c_1(P) = \dfrac{1}{k!} \sum_{\substack{\emptyset = v_0 < ... < v_k =P \\ v_i \in J(P)}} \prod_{i=0}^{n-1} c_1(v_{i+1} \setminus v_i)\]

 Since $\Omega(\emptyset, t) = 1$ we can use the inclusion-exclusion principle to get :  
 \begin{eqnarray*}
      \sum_{i=1}^{k} (-1)^{k-i}\binom{k}{i}\Omega(P;it)  &=& \sum_{\substack{\emptyset = v_0 < ... < v_k =P \\ v_i \in J(P)}} \prod_{i=0}^{k-1}\Omega(v_{i+1} \setminus v_i; t) \\
       \left[\sum_{i=1}^{k} (-1)^{k-i}\binom{k}{i}\Omega(P;it) \right]_{t^k} &=& \left[\sum_{\substack{\emptyset = v_0 < ... < v_k =P \\ v_i \in J(P)}} \prod_{i=0}^{k-1}\Omega(v_{i+1} \setminus v_i; t) \right]_{t^k} \\ 
       \sum_{i=1}^{k} (-1)^{k-i}\binom{k}{i} i^kc_{k}(P) &=& \sum_{\substack{\emptyset = v_0 < ... < v_k =P \\ v_i \in J(P)}} \prod_{i=0}^{k-1}c_1(v_{i+1} \setminus v_i; t) \\
       (k!)c_k(P) &=& \sum_{\substack{\emptyset = v_0 < ... < v_k =P \\ v_i \in J(P)}} \prod_{i=0}^{k-1}c_1(v_{i+1} \setminus v_i; t) \\
       c_k(P) &=& \frac{1}{k!}\sum_{\substack{\emptyset = v_0 < ... < v_k =P \\ v_i \in J(P)}} \prod_{i=0}^{k-1}c_1(v_{i+1} \setminus v_i; t)
 \end{eqnarray*}
\end{proof}
 
\begin{corollary}\label{corrollary:1}
	Given $P$ a poset of size $n$, recall that $a_k(P) := (|P|!)\; c_k(P)$ then for all $k \in \mathbb{N}$ : 
	\[a_k(P) = \dfrac{1}{k!}\sum_{\substack{\mu \in \mathsf{SDecompo}_{k}(P) \\  \mu = (P1,P_2,...,P_k)}} \binom{n}{\mu}e(\mu) \Pi_{i = 1}^k a_1(P_i) .\]
\end{corollary}
\begin{remark}
A nicer way to formulate Theorem~\ref{thm:1} is by not regrouping all linear extension of a given standard decomposition. Let $\text{LSDecompo}_k(P)$ the set of linear extension on standard decomposition into $k$ subposet, i.e. the set $\{(\mu, l) | \mu \in \text{SDecompo}_k(P) \text{ and } l \in \mathcal{L}_\mu\}$. As mentioned before, this set is in bijection with the set of increasing sequences of ideals of $P$ of size $k$ hence: 

    \[c_k(P) = \dfrac{1}{k!}\sum_{\substack{(\mu, l)  \in \text{LSDecompo}_{k}(P) \\  \mu = \{P_1,P_2,...,P_k\}}}\Pi_{i = 1}^k c_1(P_i) .\]
\end{remark}

\section{The weighted-sum formula}

The goal of this section is to prove the weighted-sum formula (Theorem~\ref{thm:3}). This formula expresses the coefficient $\left[\Omega(Z_n; t-1/2)\right]_{t^k}$ as a weighted sum over all the subposet decomposition of $Z_n$ into $k$ subposet of odd size. The weight can be computed by looking at the Taylor expansion of the series defined in Table~\ref{master table}. The strategy of the proof is the following. In Section \ref{subsect:constant_linear}, using the Kreweras formula, we compute the linear coefficient and the constant coefficient of $\Omega(Z_n; t-1/2)$. In Section~\ref{sec:43}, we express all terms of the sum of $ \left[\Omega(Z_n;t-1/2)\right]_{t^k} = \big[\sum_{i \geq k} c_i(Z_n)(t-1/2)^i\big]_{t^k}$ by defining a map from the decompositions into $i$ subposets to the decompositions into $k$ subposets where $i\geq k$. This leads to the definition of a weight function $w$ on the subposets. In Section \ref{elementarydecompo}, we define the notion of  elementary decomposition which is the cornerstone of all the decompositions, then we show that for elementary ones, the weight $w$ can be computed using the results of Section \ref{subsect:constant_linear}. Finally, in Section \ref{multiplicative_behavior} we show how we can split every decomposition into elementary decompositions and show that in this case, the weight of every decomposition is multiplicative under this split. We write \new{$f_{n,k}$} as the coefficients of $\Omega(Z_n;t-1/2)$. More precisely, $\Omega(Z_n;t-1/2) = \sum_{k=1}^n f_{n,k}t^k$. 

\subsection{Generating function}
From here on, we will use several generating functions and their corresponding Taylor expansion coefficients at $z=0$. Table~\ref{master table} gives the correspondence between the coefficients and the ordinary generating function.

\begin{table}[h]
\renewcommand{\arraystretch}{2}
\begin{tabular}{cccc}
\hline
Name  & Closed formula  & Coefficient formula \\
\hline
$W_{+}(z) = \sum_{k \geq 0} w_{+}(2k+1)z^{2k+1}$  & ${2\arcsin\left(z/2\right)}/{\sqrt{1 - z^2/4}}$ &    $w_{+}(2k+1) = \dfrac{(k!)^2}{(2k+1)!}$  \\
$W_{-}(z) = \sum_{k \geq 0} w_{-}(2k+1)z^{2k+1}$ & $2\arcsin(z/2) \sqrt{1 - z^2/4}$  &   $w_{-}(2k+1) = \begin{cases} 1& k=0  \\
\dfrac{-(k!)^2}{(2k)(2k+1)!} &k\geq 1\end{cases}$ \\
$C(z) = \sum_{k \geq1} C_{2k}z^{2k}$  & $4\arcsin^2(z/2)$  &  $C_{2k} = \frac{2((k-1)!)^2}{(2k)!}$  \\
$R(z) = \sum_{k \geq 0} R_{2k}z^{2k}$ & $\sqrt{1-z^2/4}$  &   $R_{2n}= \begin{cases} 1& n=0  \\  R_{2n} = \frac{-1}{16^n} \times \frac{2}{n}  \binom{2n-2}{n-1} & n\geq 1 \end{cases}$ \\
 $U(z) = \sum_{k \geq 0} U_{2k+1}z^{2k+1}$ & $2 \arcsin(\frac{z}{2})$ &  $U_{2k+1} = \frac{1}{(2k+1)2^{^{4k-1}}}\binom{2k}{k}$ \\
\hline
\end{tabular}
\caption{Table of the series occurring in the proof of the weighted sum formula}
\label{master table}

\end{table}

\begin{remark}
     Notice that $W_{+}, C, U$ have non-negative coefficients, whereas $W_{-}$ and $R$ have non-positive coefficients (except for $w_{-}(1)$ and $R_{0}$, respectively)
\end{remark}

\subsection{Constant and linear coefficients}
\label{subsect:constant_linear}

As mentioned in  \cite{ferroni2025skew}, the order polynomial $\Omega_{\lambda / \mu}(t)$ of the cell poset of a  skew shape $\lambda / \mu$ counts the number of plane partition of shape $\lambda / \mu$ with entry less than equal to $t-1$. 

\begin{theorem}[{Kreweras \cite{kreweras1965classe}}]
    Let $\lambda / \mu$ a skew shape, then the order polynomial is : 
    \[\Omega_{\lambda / \mu}(t) = P_{\lambda / \mu}(t-1) = \det \left[\binom{\lambda_i - \mu_j + t -1}{\lambda_i - \mu_j - i + j} \right]_{i,j=1}^l.\]
    
\end{theorem}
Zigzag posets are skew-shape posets : 
\begin{enumerate}
    \item if $n$ is odd, let $\lambda_n = (\frac{n+1}{2}, \frac{n-1}{2}, \cdots, 1)$ and $\mu_n = (\frac{n-1}{2}, \frac{n-3}{2}, \cdots, 1, 0, 0 )$ then $\Omega_n(t) = \Omega_{\lambda_n / \mu_n}(t)$
    \item if $n$ is even, let $\lambda_n = (\frac{n}{2}, \frac{n-2}{2}, \cdots, 0)$ and $\mu_n = (\frac{n-2}{2}, \frac{n-4}{2}, \cdots, 1, 0, 0, 0 )$ then $\Omega_n(t) = \Omega_{\lambda_n / \mu_n}(t)$
   
\end{enumerate}

Let $\textcolor{blue}{A(n)}=(A_{i,j}(n))_{i,j=1}^{\lceil n/2\rceil}$ be the matrix in the Kreweras formula associated to the zigzag $Z_n$. We collect elementary properties of this matrix.

\begin{proposition}
For the matrix $A(n)=(A_{i,j}(n))_{i,j=1}^{\lceil n/2\rceil}$ defined above : 
\begin{enumerate}[(i)]
    \item if $i >  j+1, \quad \binom{\lambda_i - \mu_j + t}{\lambda_i - \mu_j - i + j} = 0 $, \quad if $i = j+1 \quad \binom{\lambda_i - \mu_j + t}{\lambda_i - \mu_j - i + j} = 1 $ 
     \item if $n$ is even, $A_{1,j}(n)= \binom{t+j}{2j}$
    \item if $n$ is odd, $A_{1,\frac{n}{2}}(d)= \binom{t+ \frac{n-1}{2}}{d}$ and for $j < \frac{n}{2}$,  $A_{1,j}(n)= \binom{t+j}{2j}$ 
    \item $A(n-2)$ is the minor of $A(n)$ when the first row and column have been removed
\end{enumerate}
\end{proposition}

Using these properties, we obtain the following minor decomposition.

\begin{proposition}\label{lemma:8}
    If $k \geq 1$, then $\Omega(Z_{2k+1};t) = (-1)^k\binom{t+k}{2k+1} + \sum_{i=1}^k (-1)^{i+1}\binom{t+i}{2i} \Omega(Z_{2(k-i)+1};t)$
and     $\Omega(Z_{2k};t) = \sum_{i=1}^k (-1)^{i+1}\binom{t+i}{2i} \Omega(Z_{2(k-i)};t).$
\end{proposition} 

\vspace{6mm}
We are now able to compute the constant and linear coefficient of $\Omega(Z_n;t-1/2)$. The technical proof of these results can be found in Appendix.


\begin{restatable}{lemma}{oddlinear}
\label{lemma:10}
 Let $n \geq 1$. then \[f_{(2n+1),0} = 0 \quad  \text{ and } \quad f_{(2n+1),1} = w_{-}(2n+1).\]
\end{restatable}
    

\begin{restatable}{lemma}{evenconstant}
\label{lemma:11}

    Let $n \geq 1$ : then 
    \[f_{2n,0} = R_{2n}  \quad \text{  and  } \quad f_{2n,1} = 0\]
\end{restatable}
\subsection{Homogenization, decomposition into elementary blocks}
\label{sec:43}
We continue calculating the coefficients $f_{n,k}$. By linearity of extracting coefficients, we have that 
\begin{align*}
f_{n,k} &= n!\,\left[\Omega(Z_n; t-1/2)\right]_{t^k}
= n!\,\left[\sum_{i=0}^{n-k} c_{k+i}(Z_n) (t-1/2)^i\right]_{t^k}\\
&= n!\sum_{i=0}^{n-k} c_{k+i}(Z_n)\binom{k+i}{i}(-2)^{i-k}. 
\end{align*}
We define a function \new{$T$} which will allow us to make a link between decomposition into $k+i$ subposet and decomposition into $k$ subposet, therefore allowing us to compute the sum. 
Define the function $T$ : 
\begin{align*}
T \colon  \quad \mathsf{LSDecompo}_k(Z_n) &\longrightarrow\mathsf{LSDecompo}_{k-1}(Z_n) \\
               (\mu_1,l_1) &\longmapsto(\mu_2, l_2)
\end{align*}
with the convention $\text{LSDecompo}_{0}(Z_n)  = (\emptyset, \emptyset) $ \\
as: \\
If $(\mu_1, l_1) \in \text{LSDecompo}_1(Z_n)$, then $T((\mu_1, l_1)) = (\emptyset, \emptyset)$ \\
Let $x$ be the subposet of $(\mu_1,l_1)$ with highest label. Let $y$ be the subposet of $(\mu_1,l_1)$ dominated by $x$ with highest label. then let $T(\mu_1,l_1) $ be the result of the merge between $x$ and $y$. \\

Also, we define $T^m = T \circ T \cdots \circ T $, the m'th iterate of $T$. 
\begin{figure}[h]
    \centering
    \includegraphics[width=0.5\linewidth]{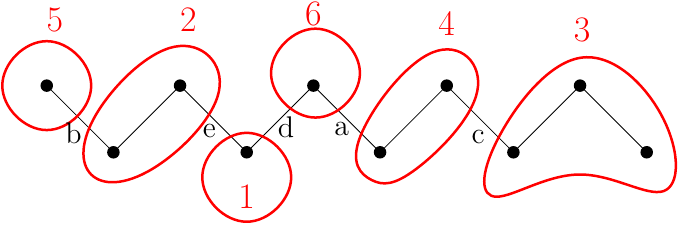}
    \caption{Illustration of the function $T$. In numerals, the ordered standard decomposition. And with letters, the order of the merges (alphabetical order). For example, if you apply $T$ twice, you will obtain $(\{5\}, \{2\}, \{ 1\}, \{6, 4\} , \{3\})$  and then $(\{5, 2\}, \{ 1\}, \{6, 4\} , \{3\})$. }
\end{figure}

\begin{remark}
    In order to reverse $T^m$, we only need to store the order of the merges. 
\end{remark}

\vspace{10mm}
Let $(\mu, l) \in \mathsf{LSDecompo}_{k}(Z_n)$. 
Let \[\textcolor{blue}{s(n,i, \mu, l)} := \displaystyle \sum_{\substack{(\mu_1, l_1) \in \mathsf{LSDecompo}_k(Z_n) \\ T^i(\mu_1,l_1) = (\mu, l)}} 2^{-i}\prod_{P_j \in \mu_1}c_1(P_j). \]

\begin{lemma}
    For all $l, l' \in \mathcal{L}_\mu$, $s(n,i, \mu, l) = s(n,i, \mu, l') $
\end{lemma}
\begin{proof}

Let $\phi : T^{-i}(\mu,l) \longrightarrow T^{-i}(\mu,l')$
defined as following. 
\begin{enumerate}
    \item Start with an element $(\mu_1,l_1) \in T^{-i}(\mu,l)$
    \item Apply $T^i$ on $(\mu_1,l_1)$ and keep  stored the information of the order of the merges. 
    \item You now have $(\mu,l)$. Replace it by $(\mu,l')$
    \item Use the reversed order to rewind the algorithm and get an element $(\mu_1,l'_1) \in T^{-i}(\mu,l')$
\end{enumerate}
$\phi$ is a bijection that preserve the decomposition. Therefore : 
\begin{eqnarray*}
    s(n,i, \mu, l)  &=& s(n,i,\mu, l')
\end{eqnarray*} 
Hence, we define $\textcolor{blue}{s(n,i,\mu)} = s(n,i,\mu, l)$ for any $l \in \mathcal{L}_{\mu}$.
\end{proof}

Also, let $(\mu_0, l_0)$ be the only standard decomposition into one subposet (with the trivial linear extension). \\
We have $T^{-(k-1)}(\mu_0,l_0) = T^{-k}(\emptyset, \emptyset)$
In particular, when $i=n$, we get : 
\begin{lemma}\label{lemma:14}
    \[f_{n,0} = \sum_{i=0}^n \dfrac{(-1)^i}{i!} s(n,i, \emptyset, \emptyset) = 2\sum_{i=0}^n \dfrac{(-1)^i}{i!} s(n,(i-1), \mu_0, l_0). \]
\end{lemma}

\begin{definition}
For $(\mu, l) \in \mathsf{LSDecompo}_k(Z_n)$, let : 
\begin{align}
\textcolor{blue}{w(\mu,l)} &:=\sum_{i=0}^{n-k} \dfrac{(-1)^i}{i!} s(n,i, \mu,l) \label{eq:def:Hmul}.\\ 
\textcolor{blue}{w(\mu)} &:=\sum_{i=0}^{n-k} \dfrac{(-1)^i}{i!} s(n,i, \mu). \label{eq:def:Hmu}
\end{align}
\end{definition}
We can therefore write $f_{n,k}$ as in the following lemma.  

\begin{lemma}\label{lemma:13}
    \begin{equation}
        f_{n,k} = \dfrac{n!}{k!} \sum_{(\mu,l) \in \mathsf{LSDecompo}_{k}(Z_n)}w(\mu,l) = \dfrac{n!}{k!} \sum_{\mu \in \mathsf{LSDecompo}_{k}(Z_n)}w(\mu)e(\mu). \
    \end{equation}
\end{lemma}
\begin{proof}

\begin{eqnarray*} f_{n,k} &=& n!\sum_{i=0}^{\,n-k} \binom{k+i}{\,i}\,(-2)^{-i}\, c_{\,k+i,\,n} \;=\; n!\sum_{i=0}^{\,n-k} \binom{k+i}{\,i}(-2)^{-i} \frac{1}{(k+i)!} \sum_{(\mu_1,l_1)\in \text{LSDecompo}(Z_n)} \prod_{P_j\in\mu_1} c_1(P_j) \\[1em] &=& n!\sum_{i=0}^{\,n-k} \sum_{(\mu,l)\in \text{LSDecompo}_{k}(Z_n)} (-1)^i \binom{k+i}{k}\, \frac{1}{(k+i)!}\, s(n,i,\mu,l) \\[1em] &=& \frac{n!}{k!} \sum_{i=0}^{\,n-k} \sum_{(\mu,l)\in \text{LSDecompo}_{k}(Z_n)} \frac{(-1)^i}{i!}\, s(n,i,\mu,l) \\[1em] &=& \frac{n!}{k!} \sum_{(\mu,l)\in \text{LSDecompo}_{k}(Z_n)} w(\mu,l) \;=\; \frac{n!}{k!} \sum_{\mu\in\text{SDecompo}_{k}(Z_n)} w(\mu)\,e(\mu). \end{eqnarray*}
\end{proof}

\subsection{Computation of the weights \texorpdfstring{$w(\mu)$}{w(mu)} for elementary decompositions}
\label{elementarydecompo}
An elementary decomposition is a decomposition in which all but one subposets have size $1$. We show in this section that if $\mu$ is an elementary decomposition, $w(\mu)$ is the linear term or the constant term of an $\Omega(Z_r;t-1/2)$, therefore, we know its value. 

\begin{definition}
    Let $n$ a positive integer, $m \leq n$ a positive integer, and let $\mu$ a decomposition of $Z_n$. We say that $\mu$ is an \new{elementary decomposition of size $m$} if  $\mu$ is composed of $n-m$ subposets of size one and one subposet of size $m$.
\end{definition}
\begin{figure}[h]
    \centering
    \includegraphics[width=0.5\linewidth]{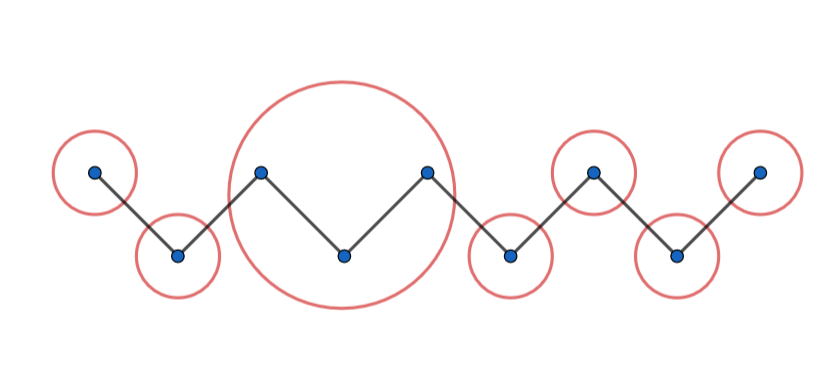}
    \caption{Example of an elementary decomposition of size $3$. All the blocks except one have size one.}
    \label{fig:size3}
\end{figure}

\begin{definition}
    An elementary decomposition $\mu$ of size $m$ is \new{special} if and only if $n$ is even, and the $m$-size subposet of $\mu$ is the rightmost subposet in $\mu$. 
\end{definition}

\vspace{5mm}
Let us try to compute $w(\mu)$ for elementary decompositions. We will then show that $w(\mu)$ has a multiplicative behaviour thus, we will be able to calculate $w(\mu)$ by decomposing $\mu$ into elementary decompositions. 

\begin{restatable}{lemma}{Hformule}
\label{lemma:16}
    Let $n$ a positive integer, let $m \leq n$ a positive integer and let $\mu$ be an elementary decomposition of $Z_n$ of size $m$. Then for $w(\mu)$ defined in \eqref{eq:def:Hmu}:
    \begin{enumerate}
        \item If $m$ is even :
        \[ w(\mu)  = 0. \]
        \item If $m$ is odd and the size-$m$ subposet has more maximal vertices than minimal vertices : 
        \[w(\mu)   = w_{-}(m).\]
        \item If $m$ is odd and the size-$m$ subposet has more minimal vertices than maximal vertices  and $\mu$ {\bf is  not special}:
        \[w(\mu)   = w_{+}(m). \]
        \item If $m$ is odd and the size-$m$ subposet has more minimal vertices than maximal vertices  and $\mu$ {\bf is special }:
        \[w(\mu)   = U_{m}.\]
    \end{enumerate}
\end{restatable}
\begin{proof}
    (see Appendix)
\end{proof}
\vspace{10mm}

\subsection{Multiplicative behavior}
\label{multiplicative_behavior}
Let $\mu_1, \mu_2$ two decomposition of $Z_n$ such that for each vertex $v$ of $Z_n$, $v$ is in at most one subposet of size greater or equal of 2 (either in  $\mu_1$ or in $\mu_2$). Let $\mu = M(\mu_1, \mu_2)$ the fusion of those two decompositions, i.e, for each element $v$, $v$ is in the subposet   (either in  $\mu_1$ or in $\mu_2$) of maximal size. 

\begin{lemma}\label{lemma:17}Let $n \geq 1$ and $0 \leq i \leq n$. 
    \[s(n,i, \mu) = \sum_{r=0}^i \binom{i}{r} s(n,r, \mu_1) s(n,i-r, \mu_2).\]
\end{lemma}
\begin{proof}
    Decompose, according to the cardinality of $r$, the number of merges in non-singleton subposets of $\mu_1$. 
\end{proof}

The next result shows that the $w(\cdot)$ is multiplicative with the block operation $M(\cdot,\cdot)$.

\begin{corollary}\label{corro:3} 
     \[w(M(\mu_1,\mu_2)) = w(\mu_1)\cdot w(\mu_2)\]
\end{corollary}
\begin{proof}
Let $k_1, k_2$ be respectively the number of subposet in $\mu_1$ and $\mu_2$
    \begin{eqnarray*}
        w(\mu) &=& \left( \sum_{i=0}^{n-k} \dfrac{(-1)^i}{i!} s(n,i, \mu)\right) \\
        &=& \left( \sum_{i=0}^{n-k} \dfrac{(-1)^i}{i!} \sum_{r=0}^i \binom{i}{r} s(n,r, \mu_1) s(n,i-r, \mu_2)\right) \\
        &=& \left( \sum_{i=0}^{n-k_{1}} \dfrac{(-1)^i}{i!} s(n,i, \mu)\right)\left( \sum_{i=0}^{n-k_2} \dfrac{(-1)^i}{i!} s(n,i, \mu)\right) \\
        &=& w(\mu_1)\cdot w(\mu_2).
    \end{eqnarray*}
\end{proof}

\subsection{The weighted sum formula}\label{sec:55}~\\

We recall in one place the definition of the weight function for decomposition into subposets. This is a summary of of Lemma~\ref{lemma:16} and Corollary~\ref{corro:3}. 
\begin{itemize}
	\item For each subposets $P$ of even size, $w(P) = 0$.
	\item For each subposets $P$ of size $2k+1$, with more upper elements than lower elements, $w(P) = w_{-}(2k+1)$.
    \item For each {\bf special} subposets $P$ of size $2k+1$, with more lower elements than upper elements, $w(P) = U_{2k+1}$.
	\item For each {\bf non-special } subposets $P$ of size $2k+1$, with more lower elements than upper elements, $w(P) = w_{+}(2k+1)$.
\end{itemize}
For a decomposition $\mu$ into subposets, we let $w(\mu) = \displaystyle \prod_{P \in \mu} w(P)$. \\
The weight function $w$ describes the coefficients of $\Omega(Z_n; t-1/2)$. The following theorem is the main result of this paper. 
\begin{theorem}\label{thm:3}
    Let $n \geq 2$ and $1 \leq k \leq n$, then the coefficient $f_{n,k} = \left[\Omega(Z_n;t-1/2)\right]_{t^k}$  equals: 
    \begin{equation}
        f_{n,k} = \dfrac{n!}{k!}E_{k}\sum_{\mu \in \mathsf{SDecompo}_{k}(Z_n)} w(\mu).
    \end{equation}
\end{theorem}

\begin{proof}
    Using Corollary~\ref{corro:3}, lemma~\ref{lemma:16} and Lemma~\ref{lemma:13}, we have : 
    \[f_{n,k} = \frac{n!}{k!} \sum_{\mu \in \mathsf{SDecompo}_k(Z_n)} w(\mu)e(\mu)\]

    If $\mu$ has at least an even-size subposet, $w(\mu) = 0$. And if $\mu$ has only odd-size subosets $e(\mu) =e({Z_k}) = E_k$, Hence : 

    \[f_{n,k} = \dfrac{n!}{k!}E_{k}\sum_{\mu \in \mathsf{SDecompo}_{k}(Z_n)} w(\mu).\]
\end{proof}
\begin{example}
    Figure~\ref{fig:example f6_2} shows all the decomposition of $Z_6$ into two subposets, therefore $f_{6,2} = \frac{6!}{2!}E_2 (w_-(3)w_{+}(3) + w_-(5) + U_5)$.
    \begin{figure}[h]
    \centering
    \includegraphics[width=0.5\linewidth]{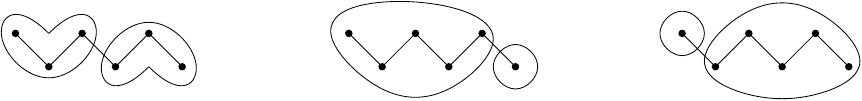}
    \caption{The $3$ ways to decompose $Z_6$ into $2$ subposets. The decomposition have respective weights $w_-(3)w_{+}(3)$, $w_-(5)$ and $U_5$.}
    \label{fig:example f6_2}
\end{figure}
\end{example}

The next two observations of Theorem~\ref{thm:3} give non-trivial properties about the coefficients $f_{n,k}$, namely a closed formula for $f_{n,n-2}$ independent of the parity of $n$, and the fact that every other coefficient of $\Omega(Z_n;t-1/2)$ is zero.

\begin{example}
Here, we show as an example how to use the weighted sum formula to compute $f_{n,n-2}$. Notice that if any decomposition of $Z_n$ into $n-2$ convex subposets of odd size is constituted as a convex subposet of size $3$ and the other subposet of size $1$. Notice that $w_{-}(3) = \frac{-1}{12}$, $w_{+}(3) = \frac{1}{6}$ and $U_3 = \frac{1}{24}$. If $n$ is odd, there are $(n-1)/2$ decompositions of weight $w_{-}(3)$ and $(n-3)/2$ decompositions of weight $w_{+}(3)$ therefore, $f_{n, n-2} = n(n-1)E_{n-2}\cdot(\frac{n-4}{12} - \frac{n-1}{24}) = \frac{n(n-1)(n-5)E_{n-2}}{24}$. If $n$ is even, there are $(n-4)/2$ decompositions of weight $w_{+}(3)$, $\frac{n-2}{2}$ decompositions of weight $w_{-}(3)$ and one decomposition of weight $U_3$. Therefore $$f_{n,n-2} = n(n-1)E_n \cdot \left(\frac{n-4}{12} - \frac{n-2}{24} + \frac{1}{24}\right) = n(n-1)(n-5)E_n/{24}.$$ At first glance, it is quite surprising that the formula is the same for $n$ even or odd.
\end{example}

\begin{remark}
  If $j$ is an odd integer, any element of $\mathsf{SDecompo}_{n-j}(Z_n)$ has a subposet of even size. Thus, $[\;\Omega(Z_n;t-1/2)]_{t^{n-j}}= 0$. So, if $n$ is odd or even, $\Omega(Z_n;t-1/2)$ is an odd or even polynomial, respectively. 
\end{remark}

\section{Analysis of the generating function of the weighted sum}

We want to be able to makes some inequalities with the coefficients from \ref{thm:3} in order to prove the Watanabe--Yoshida conjecture. First, let us write \[\textcolor{blue}{g_{n,k}} := \dfrac{k!f_{n,k}}{n!E_k} = \sum_{\mu \in \mathsf{SDecompo}_{k}(Z_n)} w(\mu).\] We write $\textcolor{blue}{G_n(x)} = \sum_{k \geq 0} g_{n,k} x^k$ and $\textcolor{blue}{G(z;x)} := \sum_{n \geq 0} G_n(x) z^{n}$.
We define $g_{0,0} := 1$
\begin{lemma}\label{lemma:18}
    Let $n$ a positive integer and $k \geq 3$ 
        \[g_{n,k} = \sum_{\substack{i=2 \\  i \text{ even }}}^nC_i \cdot g_{n-i,k-2}.\]
Therefore : 
    \begin{eqnarray*}
    G_{\text{odd}}(z;x) = x\dfrac{W_{-}(z)}{1-x^2C(z)}, &\quad \quad&G_{\text{even}}(z;x) = \dfrac{R(z)}{1-x^2C(z)},
\end{eqnarray*}
and \[G(z;x) =   \sqrt{1-z^2/4} \cdot\left( 1 + 2x\arcsin (z/2)\right)/(1-4x^2 \arcsin(z/2))\]
\end{lemma}
\begin{proof}
    Let $i$ be the sum of the size of the two leftmost subposets.
    Therefore, we need to evaluate : 
    \[\sum_{\substack{j=1 \\ j\text{ odd}}}^i w_{-}(j)w_{+}(i-j) = \left[W_{-}(z)W_{+}(z)\right]_{z^i} = [4\arcsin(z/2)]_{z^i} = C_i.\]

For $G_{\text{odd}}$, we see that $[G(z;x)]_{x^1}$ is $W_{-}(z)$ and use Lemma~\ref{lemma:18}. For $G_{\text{even}}$ using that $[G(z;x)]_{x^0} = R(z)$ (Lemma~\ref{lemma:11}). 
And calculate that $[G(z;x)]_{x^2} = W_+(z)U(z) = R(z)C(z)$, we have 
$G_{\text{even}}(z;x) = R(z) + \frac{x^2R(z)C(z)}{1-x^2C(z)} = \frac{R(z)}{1-x^2C(z)}$
\end{proof}

Using the previous result, Theorem~\ref{thm:3} can be restated as a Hadamard product $\odot$. Given two series $A(t) = \sum_{n\geq0} a_nz^t$ and $B(t)= \sum_{n \geq 0} b_n t^n$. The Hadamard product is defined as $(A\odot B )(t) = \sum_{n \geq0} a_nb_nt^n$.  

\begin{theorem}\label{thm:hadamardprod}
    For $d\geq1$ and for all $t$ : 
    \begin{equation}
        \sum_{n \geq 0} \;\Omega(Z_n;t-1/2)\frac{z^n}{n!} = \left( \tan(t) + \sec(t) \right) \odot G(z;t).
    \end{equation}
\end{theorem}
\begin{proof}
    \begin{eqnarray*}
        \sum_{n \geq 0} \;\Omega(Z_n;t-1/2)\frac{z^n}{n!} &=& \sum_{n \geq 0 }\sum_{k =0}^n f_{n,k}t^k\frac{z^n}{n!} \\
        &=&\sum_{n \geq 0 }z^n\sum_{k =0}^n g_{n,k}\frac{E_kt^k}{k!} \\
        &=& \sum_{n \geq0}z^n \left( G_n(t) \odot (\tan(t) + sec(t)\right)\\
        &=& G(z;t) \odot (\tan(t)+\sec(t))
    \end{eqnarray*}
\end{proof}

\begin{example}
    For example, to compute $\Omega(Z_6;t-1/2)$, we have 
    \begin{align*}
    G_{6}(t) &= &\frac{-1}{1024}&&+\frac{-41}{5760}t^2 &&+ \frac{1}{24}t^4 &&+ t^6,\\
    \tan(t) + \sec(t) &=  &1&+ \quad t &+ \frac{1}{2}t^2 &+ \quad \frac{1}{3}t^3 &+ \frac{5}{24}t^4 &+ \quad \frac{2}{15}t^5 &+ \frac{61}{720}t^6 &+ O(t^7),\\
    \Omega(Z_6;t-1/2) &= &\frac{-1}{1024} &&+ \frac{-41}{5760 \cdot 2}t^2  &&+ \frac{5}{24 \cdot 24}t^4 &&+ \frac{61}{720}t^6.
    \end{align*}
    Lastly, shifting the variable $t$ we have that $\Omega(Z_6;t)=\frac{1}{6!}(12t + 64t^2 + 165t^3 + 235t^4 + 183t^5 + 61t^6)$.
\end{example}

We start by recalling the well-known asymptotic formula for the Euler numbers.
\subsection{Asymptotic analysis}

\begin{lemma}[e.g \cite{stanley2009survey}]\label{lemma:19}
For all $n\geq 1$ we have 
    \[ \frac{E_n}{n!} = 2 \cdot \epsilon(n)\cdot \left(\dfrac{2}{\pi}\right)^{n+1}, \]    with $\epsilon(n) = \sum_{k\geq0} \dfrac{(-1)^{k(n+1)}} {(2k+1)^{n+1}}$.
\end{lemma}
      
\begin{definition}
    let $x \in \mathbb{C}$, and $f$ a function holomorphic at $z=0$, then the following function is well-defined :  
    
    \begin{align*}
        \textcolor{blue}{\mathsf{Alt}_x} :\mathcal{F} \rightarrow \mathbb{C}, \qquad 
    &f \mapsto \sum_{k \geq 0 }\frac{(-1)^k}{(2k+1)}f\left(\frac{x(-1)^k}{(2k+1)}\right).
    \end{align*}
And 
   \[\textcolor{blue}{S(n)} := \mathsf{Alt}_x(x^n)_{x=1} = \sum_{k \geq0}\frac{(-1)^{k(n+1)}}{(2k+1)^{(n+1)}}.\]
\end{definition}

Using this notation and Lemma~\ref{lemma:16}, we can write  ${E_n}/{n!} = [2\cdot \mathsf{Alt}_x(x^n)]_{x=\frac{2}{\pi}}$. Now, let us expand $\Omega(Z_n;t-\frac{1}{2})$.  Let 
  $\textcolor{blue}{H_n(x)} := \sum_{k=0}^{n-2}g_{n,k}x^k = G_{n}(x) - x^n$.

\begin{proposition}
    \begin{equation} \label{eq:Omega in terms of Alt}
        n!\cdot \,\Omega_n\!\left(t-\dfrac{1}{2}\right)=E_n t^{\,n}
        + \dfrac{4\, n!}{\pi}
           \left(\mathsf{Alt}_t\!\left(
               H_n\!\left(\dfrac{2t}{\pi}\right)
           \right) \right),
    \end{equation}

\end{proposition}
\begin{proof}
We use Lemma~\ref{lemma:19} to rewrite ${E_k}/{k!}$ with the $\mathsf{Alt}$ notation. 
\begin{eqnarray*}
    n! \cdot \,\Omega_n\!\left(t-\dfrac{1}{2}\right)
    &=& E_n t^{\,n}
        + \sum_{k=0}^{\,n-2} f_{n,k}\, t^{k} \\
    &=& E_n t^{\,n}
        + n! \dfrac{4}{\pi}
           \left( \sum_{k=0}^{\,n-2}
           \left(\dfrac{2}{\pi}\right)^{k}
           \epsilon(k)\, t^{k}\, g_{n,k}(n,x)
           \right) \\
    &=& E_n t^{\,n}
        + \dfrac{4\, n!}{\pi}
           \left( \sum_{i \geq 0}
           \dfrac{(-1)^i}{2i+1}
           \sum_{k=0}^{\,n-2}
           \dfrac{2^{k} (-1)^kt^k}{\left(\pi(2i+1)\right)^{k}}
           g_{n,k}
           \right) \\
    &=& E_n t^{\,n}
        + \dfrac{4\, n!}{\pi}
           \left(\mathsf{Alt}_t\!\left(
               H_n\!\left(\dfrac{2t}{\pi}\right)
           \right) \right),
\end{eqnarray*}
as desired.
\end{proof}

The last term of the right-hand-side of \eqref{eq:Omega in terms of Alt}  is positive for $t$ and $n$ big enough. We postpone the proof of this fact to Section~\ref{sec6}.

\begin{lemma}\label{prop:11}
    For $n \geq 8, \; t \geq \frac{3}{2}$ : 
    \[\mathsf{Alt}_t\left(H_n\left(\frac{2t}{\pi}\right) \right)\ > 0. \]
\end{lemma}
\begin{proof}
See Section~\ref{sec6}.

\end{proof}

We use the inequality above to give a lower bound of the shifted order polynomial of zigzag posets. 

\begin{corollary}\label{corro:finalinequality}
    For $n \geq 8$ and $p \geq 3$, 
    \[\Omega(Z_n;p/2-1/2) > \frac{E_n}{n!}\cdot \left(\frac{p}{2} \right)^n.\]
\end{corollary}

We compute the multiplicity $e_{HK}(\mathcal{A}_{p,n})$ using \eqref{eq:pakformula} and check computationally  using {\tt SageMath} \cite{sagemath} that indeed $e_{HK}(\mathcal{A}_{p,n}) \geq 1+E_n/n!$ in those cases. This finishes our proof of the third part of the Watanabe--Yoshida conjecture \cite{watanabe2005hilbert} (Conjecture~\ref{conj:1}) realizing the Ehrhart theory approach of \cite{pak2025hilbertkunzmultiplicityquadricsehrhart}.

\begin{theorem} \label{thm:proof 3rd part WS conj}
For all $p\geq 3$ and $n\geq 2$, we have that 
    \[
    e_{HK}(\mathcal{A}_{p,n})  \geq 1+\frac{E_n}{n!}.\]
\end{theorem}

\section{Proof of Lemma~\ref{prop:11}}\label{sec6}
In this section, we prove that $\mathsf{Alt}_x(H_n(x)) \ge 0$ for all $x \ge 3/\pi$ and all integers $n \ge 8$. The argument relies on the decomposition of $H_n(x)$ as a sum of functions with positive expansions and a careful analysis of the remaining terms.

\subsection{Simplification of \texorpdfstring{$H_n$}{H\_n}}

Let 
\[
H(z;x) := \sum_{n \ge 0} H_n(x) z^n = \sum_{n \ge 0} \left( G_n(x) - x^n \right) z^n.
\]
We start by writing $H(z;x)$ as a sum of two parts which are easier to analyze. The decomposition differs depending on whether $n$ is odd or even.

\begin{itemize}
\item \textbf{If $n$ is odd:} 
\begin{align*}
H_{\text{odd}}(z;x) &= \frac{x W_{-}(z)}{1-x^2 C(z)} - \frac{x z}{1-(x z)^2} \\
&= z x \left( \frac{(C(z)-z^2) + z(W_{-}(z)-z)}{(1-x^2 C(z)) z^2} + \frac{(C(z)-z^2)(2 x^2 z^2 -1)}{(1-x^2 C(z)) z^2 (1-x^2 z^2)} \right)
\end{align*}

Define 
\[
U(z;x) := \frac{(C(z)-z^2) + z(W_{-}(z)-z)}{(1-x^2 C(z)) z^2}, \quad 
B(z;x) := \frac{(C(z)-z^2)}{(1-x^2 C(z)) (1-x^2 z^2)}.
\]

Then 
\begin{equation}\label{eq:H_odd}
H_{\text{odd}}(z;x) = z x \left( U(z;x) + \frac{2x^2 z^2 -1}{z^2} B(z;x) \right).
\end{equation}

Since $w_{-}(2r+1) = \frac{-(r!)^2}{2r \cdot (2r+1)!}$ and
\[
C_{2r+2} = \frac{2((r+1)!)^2}{(r+1)^2 (2r+2)!} = \frac{(r!)^2}{(2r+1)! (r+1)} = \frac{-2r}{r+1} \cdot w(2r+1),
\]
we have $C_{2r+2} \ge - w(2r+1)$ for all $r \ge 1$. Hence, $(C(z)-z^2) + z(W_{-}(z)-z)$ has non-negative coefficients in its Taylor expansion. Therefore, $U(z;x)$ has a non-negative Taylor expansion at $z=0$, $x=0$. 

\vspace{2mm}

\item \textbf{If $n$ is even:} 

Let 
\[
H_{\text{even}}(z;x) = G_{\text{even}}(z;x) - \frac{1}{1-(x z)^2}.
\]

Then
\[
H_{\text{even}}(z;x) = \frac{\frac{3}{2}(C(z)-z^2) + z^2 (R(z)-1)}{(1-x^2 C(z)) z^2} + \frac{(C(z)-z^2)(\frac{5}{2} x^2 z^2 - \frac{3}{2})}{(1-x^2 C(z)) z^2 (1-x^2 z^2)}.
\]

Define 
\[
\hat{U}(z;x) := \frac{\frac{3}{2}(C(z)-z^2) + z^2 (R(z)-1)}{(1-x^2 C(z)) z^2}.
\]

Then
\begin{equation}\label{eq:H_even}
H_{\text{even}}(z;x) = \hat{U}(z;x) + \frac{\frac{5}{2} x^2 z^2 - \frac{3}{2}}{z^2} B(z;x).
\end{equation}

Since $R_{n+2} = \frac{n^2-1}{8 n^2} R_n$ and $C_{n+2} = \frac{n^2}{(n+2)(n+1)} C_n$, we have $\frac{C_{n+2}}{R_n} \ge 0$ for all $n \ge 2$. Therefore, the function $\hat{U}(z;x)$ has a non-negative Taylor expansion at $z=0, x=0$.

\end{itemize}

\subsection{\texorpdfstring{Analysis of $B_n(x)$}{Analysis of Bn(x)}}

\begin{lemma}\label{lemma:Bn_formula}
\[
B_n(x) = \sum_{(i_1,\ldots,i_k)\neq (1,\ldots,1)} x^{2k-2} \prod_{j=1}^k C_{2i_j},
\]
where the sum is over tuples $(i_1,\ldots,i_k)$ such that $2 i_1 + \cdots + 2 i_k = n$.
\end{lemma}

\begin{proof}
This follows from the partial fraction decomposition
\[
B(z;x) = \frac{C(z)}{1-x^2 C(z)} - \frac{z^2}{1-x^2 z^2}.
\]
\end{proof}

\begin{corollary}\label{corro:Bn_from_Bn-2}
\[
B_n(x) = x^2 B_{n-2}(x) + x^2 \sum_{\substack{i=4 \\ i \text{ even}}}^{n} C_i \left( B_{n-i}(x) + x^{\,n-2-i} \right).
\]
\end{corollary}

\subsection{\texorpdfstring{Study of $B_n$ when $n$ is odd}{Study of Bn when n is odd}}\label{sec:odd}

Recall~\eqref{eq:H_odd}. 

\begin{proposition}
For all $x \ge 0$, $\mathsf{Alt}_x(U_n(x)) \ge 0$.
\end{proposition}

\begin{proof}
The function $U_n(x)$ has a positive expansion, and for all $x \ge 0$, $\mathsf{Alt}_x(x^n) \ge 0$. Linearity of $\mathsf{Alt}_x$ completes the proof.
\end{proof}

Now consider the Taylor expansion of $(2 x^2 z^2 -1) B(z;x)/z^2$. From Corollary~\ref{corro:Bn_from_Bn-2},
\begin{align*}
B_n(x) - 2 x^2 B_{n-2}(x) &= - x^2 B_{n-2}(x) + x^2 \sum_{\substack{i=4 \\ i \text{ even}}}^{n} C_i \big( B_{n-i}(x) + x^{n-2-i} \big) \\
&= -x^2 \left( x^2 B_{n-4}(x) + x^2 \sum_{\substack{i=4 \\ i \text{ even}}}^{n-2} C_i \big( B_{n-2-i}(x) + x^{n-4-i} \big) \right) + x^2 \sum_{\substack{i=4 \\ i \text{ even}}}^{n} C_i \big( B_{n-i}(x) + x^{n-2-i} \big) \\
&= x^2 \sum_{\substack{i=0 \\ i \text{ even}}}^{n-4} B_i(x) (C_{n-i} - x^2 C_{n-2-i}) + C_n.
\end{align*}

Since $B_2(x)=0$ and $B_4(x)=1/12$, we have:

\begin{proposition}\label{prop:odd_formula}
For all $n \ge 8$,
\[
2 x^2 B_{n-2}(x) - B_n(x) = F_n + \sum_{\substack{i=4 \\ i \text{ even}}}^{n-4} B_i(x) \big( x^2 C_{n-2-i} - C_{n-i} \big),
\]
with $F_n :=\left(-C_n + \frac{x^2}{12} (x^2 C_{n-6} - C_{n-4}) \right)$ 
\end{proposition}

\begin{proposition}\label{prop:Alt_odd_sum}
For all $x \ge 3/\pi$ and $n \ge 2$,
\[
\mathsf{Alt}_x \Bigg( x \sum_{\substack{i=4 \\ i \text{ even}}}^{n-4} B_i(x) \cdot (x^2 C_{n-2-i} - C_{n-i}) \Bigg) \ge 0. \]
\end{proposition}

\begin{proof}
\begin{align*}
\mathsf{Alt}_x\Bigg(x \sum_{\substack{i=4 \\ i \text{ even}}}^{n-4} B_i(x) (x^2 C_{n-2-i} - C_{n-i}) \Bigg) 
&\ge x \sum_{\substack{i=4 \\ i \text{ even}}}^{n-4} B_i(x) (x^2 C_{n-2-i} - C_{n-i}) \\
&\quad - \frac{x}{(2k+1)^2} \sum_{\substack{i=4 \\ i \text{ even}}}^{n-4} \sum_{k \ge 1} |B_i(x_k)| |x_k^2 C_{n-2-i} - C_{n-i}| \\
&\ge x \sum_{\substack{i=4 \\ i \text{ even}}}^{n-4} B_i(x) \Big( (x^2 - \frac{\zeta_{\text{odd}}(2)}{4}) C_{n-2-i} + (\zeta_{\text{odd}}(2)-1) C_{n-i} \Big) \\
&\quad - \frac{x}{(2k+1)^2} \sum_{\substack{i=4 \\ i \text{ even}}}^{n-4} \sum_{k \ge 1} B_i(x) (x_k^2 C_{n-2-i} + C_{n-i}) \\
&\ge x \sum_{\substack{i=4 \\ i \text{ even}}}^{n-4} B_i(x) C_{n-2-i} \Big( x^2 - \frac{\zeta_{\text{odd}}(2)}{4} - (\zeta_{\text{odd}}(4)-1) \Big) \\
&\ge 0,
\end{align*}
where $x_k = \frac{(-1)^k x}{2k+1}$ and $\zeta_{\text{odd}}(n) := \sum_{k \ge 0} \frac{1}{(2k+1)^n}$, noticing that $\frac{\zeta_{\text{odd}}(2)}{4} + \zeta_{\text{odd}}(4) - 1 < 9/\pi^2 \le x^2$.
\end{proof}

\begin{proposition} For $x \geq \frac{3}{\pi}$ and $n \geq 8$,
    $\mathsf{Alt}_x(xR_n(x)) > 0$
\end{proposition}
\label{prop:term_residuel_odd}
\begin{proof}
    \begin{eqnarray*}
        \mathsf{Alt}_x(F_n) &=& \mathsf{Alt}_x\left(-C_nx + \frac{x^3}{12} (x^2 C_{n-6} - C_{n-4}) \right) \\
        &=&-S(1)C_nx - \frac{S(3)}{12}C_{n-4}x^3 + \frac{S(5)}{12}C_{n-6}x^5 \\
        &\geq& \left(\frac{16S(5)x^5}{3} + \frac{4S(3)}{3} -S(1)x\right)C_n
    \end{eqnarray*}
The latter polynomial is positive when $x >\frac{3}{\pi} $. 
\end{proof}

This finish the case $n$ odd, as we expressed $\mathsf{Alt}_x(H_n(x))$ as a sum of positive terms. 
subsection{\texorpdfstring{Study of $B_n$ when $n$ is even}{Study of Bn when n is even}}\label{sec:even}

Recall~\eqref{eq:H_even}:
\[
H_{\text{even}}(z;x) = \hat{U}(z;x) + \frac{\frac{5}{2} x^2 z^2 - \frac{3}{2}}{z^2} B(z;x).
\]

We analyze the term \((\frac{5}{2} x^2 z^2 - \frac{3}{2}) B(z;x)/z^2\).

\begin{proposition}\label{prop:even_formula}
For all \(n \ge 8\),
\[
\frac{5}{2} x^2 B_{n-2}(x) - \frac{3}{2} B_n(x) = \sum_{\substack{i=4 \\ i \text{ even}}}^{n-4} B_i(x) \big( x^2 C_{n-2-i} - C_{n-i} \big) + \hat{F}_n(x),
\]
where $\hat{F}_n = \left(
              -C_{n}
              - \frac{x^2}{12}C_{n-4}
              + \frac{x^4}{18}C_{n-6}
            \right)$
\end{proposition}

\begin{proof}
This follows from the same decomposition as in Proposition~\ref{prop:odd_formula}, using Corollary~\ref{corro:Bn_from_Bn-2} and rearranging terms. All sums of the form \(B_i(x) (x^2 C_{n-2-i} - C_{n-i})\) appear naturally, with remaining terms collected in \(R_n(x)\). The initial terms are handled separately because they may not fit the general summation pattern. The argument is identical to the odd case, just with coefficients \(\frac{5}{2}\) and \(\frac{3}{2}\).
\end{proof}

\begin{proposition}\label{prop:Alt_even_sum}
For all \(x \ge 3/\pi\) and \(n \ge 2\),
\[
\mathsf{Alt}_x \Bigg( x \sum_{\substack{i=4 \\ i \text{ even}}}^{n-4} B_i(x) \cdot (x^2 C_{n-2-i} - C_{n-i}) \Bigg) \ge 0.
\]
\end{proposition}

\begin{proof}
The proof is similar to Proposition~\ref{prop:Alt_odd_sum}
    \begin{eqnarray*}
        \mathsf{Alt}_x\left(\sum_{\substack{i=4 \\ i \text{ is even }}}^{n-4} B_i(x)\cdot\left(\frac{2}{3}x^2 C_{n-2-i} - C_{n-i}\right)\right) &=& \sum_{\substack{i=4 \\ i \text{ is even }}}^{n-4} \left(\mathsf{Alt}_x\left(\frac{2B_i(x)x^2C_{n-2-i}}{3}\right) - \mathsf{Alt}_x\left(B_i(x)C_{n-i}\right) \right)
    \end{eqnarray*}
Notice that when $f$ is an even function with positive coefficients $\mathsf{Alt}_x(f)$ is an alternating sum, therefore $\mathsf{Alt}_x(f) \geq f(x) - \frac{f(x_1)}{3}$.

Hence, 
\begin{eqnarray*}
        \mathsf{Alt}_x\left(\sum_{\substack{i=4 \\ i \text{ is even }}}^{\,n-4}
        B_i(x)\cdot\left(\frac{2}{3}x^2 C_{n-2-i} - C_{n-i}\right)\right)
        &\geq&
        \sum_{\substack{i=4 \\ i \text{ is even }}}^{\,n-4}
        \left(
            \left(\frac{2C_{n-2-i}}{3}\left(B_i(x)x^2-\frac{B_i(x_1)x_1^2}{3}\right)\right)
            - C_{n-i}\left( B_i(x) - \frac{B_{i}(x_1)}{3}\right)
        \right)
        \\
        &\geq&
        \sum_{\substack{i=4 \\ i \text{ is even }}}^{\,n-4}
        \left(
            \left(\frac{2C_{n-2-i}B_{i}(x)}{3}\left(x^2-\frac{x_1^2}{3}\right)\right)
            - C_{n-i}B_i(x)
        \right)
        \\
        &\geq&
        \sum_{\substack{i=4 \\ i \text{ is even }}}^{\,n-4}
        \left(
            \left(\frac{8x^2}{3}\left(1-\frac{1}{27}\right)-1\right)
            C_{n-i}B_i(x)
        \right)
        \\
        &>& 0 .
\end{eqnarray*}
where the last inequality holds since  $\frac{9}{4\sqrt{13}} \approx 0.624 < 2/\pi \leq x$.    
\end{proof}

Similarly as in Proposition~\ref{prop:term_residuel_odd}, we prove that $\mathsf{Alt}_x(\hat{F}_{n}(x))$. Hence, finishing the even case.

\section{Final remarks}

\subsection{The crown poset}
The crown poset $\mathcal{C}_{2n}$ is defined by taking the even zigzag $Z_{2n}$ and adding the cover relation $z_{2n} \preceq z_1$. The crown poset C2n is the poset of proper faces of a polygon with  $n$ vertices. The order polynomial of the crown poset have been studied in \cite{lundström2025orderpolytopescrownposets}.
All the proofs of this paper extend similarly to this poset, mutatis mutandis. To be more precise, the new weight function $\hat{w}$ is defined similar as in lemma~\ref{lemma:16}, except that no sub-decomposition is special. And we have to replace $E_n$ by the number of linear extension of the crown poset. Since $e(\mathcal{C}_{2n}) = nE_{2n-1}$  We therefore obtain an analog of Theorem~\ref{thm:3}.
\begin{proposition}

    Let $\Omega({\mathcal{C}_{2}};t-1/2) = \sum_{k=0}^{2n}\hat{f}_{2n,k}t^k$ : 
    \begin{equation}
        \hat{f}_{2n,k} = \dfrac{2(n!)}{(k-1)!}E_{k-1}\sum_{\mu \in \mathsf{SDecompo}_{k}(Z_{2n})} \hat{w}(\mu).
    \end{equation}
\end{proposition}

\subsection{Other shifts}
In this article we mainly focused on $\Omega(Z_n;t-1/2)$ due to its appearance in the Hilbert-Kunz multiplicity. However, for any shift and any posets, we can apply the procedure described in Section~\ref{sec:43}. It will lead to another weighted sum formula with other weights.

\section*{Acknowledgments} I thank Alejandro Morales for pointing out the conjecture in \cite{pak2025hilbertkunzmultiplicityquadricsehrhart} and for helpful comments and suggestions. I also thank Luis Ferroni, Marni Mishna, Igor Pak, Ben Shapiro, Ilya Smirnov and Ken-Ichi Yoshida for their helpful comments. The author was partially supported by the NSERC Discovery grant RGPIN-2024-06246.

\printbibliography

\section*{Appendix}



\oddlinear*
\begin{proof}
    Using Lemma~\ref{lemma:8}. We have : 
    \begin{eqnarray*}
        f_{(2k+1),1} &=& \left[(-1)^k \binom{t+k-\frac{1}{2}}{2k+1} + \sum_{i=1}^k (-1)^{i+1}\binom{t+i-\frac{1}{2}}{2i} \Omega(Z_{2(k-i)+1};t-\frac{1}{2})\right]_{t^1} 
    \end{eqnarray*}
Since $t \rightarrow \binom{t+i-\frac{1}{2}}{2i}$ is even, we just have to get the constant term if $\binom{t+i-\frac{1}{2}}{2i}$ and the linear term of $\Omega(Z_{2(k-i)+1};t-\frac{1}{2})$. 
The constant term of $\binom{t+i-\frac{1}{2}}{2i}$ is $\binom{i-\frac{1}{2}}{2i} = \dfrac{(-1)^i}{16^i} \binom{2i}{i}$. \\
The linear term of $\binom{t+k-\frac{1}{2}}{2k+1}$ is : 
\begin{eqnarray*}
    \left[ \binom{t+k-\frac{1}{2}}{2k+1} \right]_{t^1} &=& \left[ \dfrac{\mathrm{d}}{\mathrm{d}t} \binom{t+k-\frac{1}{2}}{2k+1} \right]_{t=0} \\
    &=& \binom{k+\frac{1}{2}}{2k+1} \times \left(\sum_{i=0}^{2k} \frac{1}{k-\frac{1}{2} - i}\right) \\
    &=&  \dfrac{1}{2k+1}\dfrac{(-1)^k}{16^k} \binom{2k}{k}
\end{eqnarray*}
Hence : 
\[f_{(2k+1),1} =\dfrac{(-1)^k}{(2k+1)16^k} \binom{2k}{k} +  \sum_{i=1}^k (-1)^i\binom{2i}{i} f_{1,2(k-i)} \]
Let $F(z) = \sum_{k\geq 0} f_{(2k+1),1}z^{2k}$, and 
$G(z) = \sum_{k\geq0} \dfrac{(-1)^k}{16^k} \binom{2k}{k}z^{2k}$,
and let $H(z) = \sum_{k\geq 0} \dfrac{(-1)^k}{(2k+1)16^k} \binom{2k}{k}z^{2k}$
We therefore have \[F(z) = H(z) + (G(z) -1)F(z)\]
                \[F(z) = \frac{H(z)}{G(z)}\]
A computation gives : $zF(z) = \arcsin\left(\frac{z}{2}\right) \sqrt{4 - z^2} =W_{-}(z)$

\end{proof}

\evenconstant*
\begin{proof}
    For $k \geq 1$ : 
    \begin{eqnarray*}
        f_{2k,0} &=& \left[ \sum_{i=1}^k (-1)^{i+1}\binom{t+i-\frac{1}{2}}{2i} \Omega(Z_{2(k-i)};t-\frac{1}{2})\right]_{t=0}  \\
        &=& \sum_{i=1}^k \frac{-1}{16^i} \binom{2i}{i}f_{2(k-i),0}
    \end{eqnarray*}
Let $F(z) = \sum_{k \geq0} f_{2k,0}z^{2k}$ \\
Let $G(z) = \sum_{k \geq 1} \frac{1}{16^k}\binom{2k}{k}z^{2k}$ \\
We can see that $G(z) = \dfrac{1}{\sqrt{1-\frac{z^2}{4}}}-1$
We have : 
\[F(z) = 1 - G(z)F(z)\]
\[F(z) = \dfrac{1}{1+G(z)} = \sqrt{1-\frac{z^2}{4}}\] 
\end{proof}

\Hformule*
\begin{proof}
    \begin{enumerate}
        \item~\\
        \begin{itemize}
            \item If $n$ is odd or $\mu$ is {\bf not} special. \\
            Let $l$ be an ordering of $\mu$. Let's analyze $T^{-k}(\mu,l)$.  Let $(\mu_1, l_1) \in T^{-k}(\mu,l)$, since all the merges from the algorithm are in the $m$-size subposet, the $k$ highest labels of $(\mu_1,l_1)$ are subposets of the $m$-size subposets. However, since $n$ is odd and $m$ is even, either the leftmost or the rightmost vertex of the $m$-size subposet is dominated by a vertex of $Z_n$. We call this vertex $v$.  Hence, the label of this vertex can't have one of the $k$ highest label in $(\mu_1,l_1)$. Therefore, the subposet of the $m$-size subposets containing $v$ has minimal vertex. Denote $j$ the size of the subposet of $(\mu_1, l_1)$ containing $v$. Since this subposet don't dominate any subposet of $(\mu_1, l_1)$, it is either of odd size or of size $m$. For the rest of the subposets of the $m$-size subposet, we can see those as $(\mu_2, l_2) \in \text{SDecompo}_{m-j-k}(m-j)$. 
            \begin{figure}[h]
            \centering
            \includegraphics[width=0.5\linewidth]{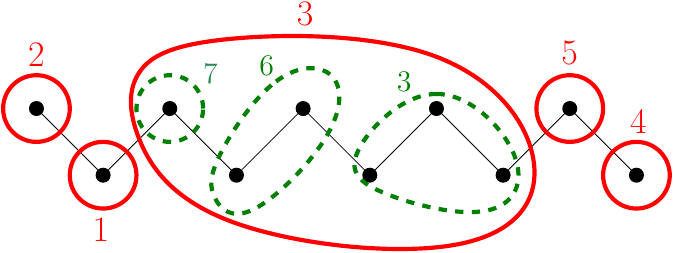}
            \caption{In red $(\mu, l)$ with $m=6$. In dashed green, an element of $T^{-3}(\mu,l)$ (the label for the subsets outside the size-$m$ subposet are the same as for $(\mu,l)$. Here, we see that $v$ is the vertex in the subposet labelled by $5$. Here $j=3$}
        \end{figure}      
        Therefore : 
        if $i \geq 1$ : 
        \begin{equation*}
        s(n,i, \mu_, l) =
            \begin{cases}
                  2 \displaystyle \sum_{\substack{j=1 \\ j \text{ odd} }}^{m-1} c_1(j) \times s(m-j, i-1,\mu_0, l_0)& \text{if } i \ge 1\\
                  c_1(m) & \text{if } i=0
            \end{cases}
        \end{equation*}
        \begin{eqnarray*}
             \sum_{i=0}^{m-1} \dfrac{(-1)^i}{i!} s(n,i,\mu, l)  &=&  c_1(m) + \sum_{i=1}^{m-1} 2\sum_{\substack{j=1 \\ j \text{ odd}}}^{m-1-i} \frac{(-1)^i}{i!} c_1(j) \times s(m-j, i-1,\mu_0, l_0) \\
             &=& c_1(m) + \sum_{\substack{j=1 \\ j \text{ odd}}}^{m-1} c_1(j)\sum_{i=0}^{m-1-j}\frac{(-1)^i}{i!}2s(m-j, i-1,\mu_0, l_0) \\
             &=& c_1(m) + \sum_{\substack{j=1 \\ j \text{ odd}}}^{m-1} c_1(j)f_{m-j,0} \\
             &=& c_1(m) + c_1(m-1)\times \frac{-1}{2} \\
             &=& 0 
        \end{eqnarray*}
        \\
        \item If $n$ is even and $\mu$ is special. In this case, neither the leftmost nor the rightmost vertex is dominated by an element outside of the $m-$size subposet. Therefore :
        \[s(n,i, \mu, l) = s(m,i, \mu_0, l_0)\]

        Hence : 
        \begin{eqnarray*}
            \sum_{i=0}^{m-1} \dfrac{(-1)^i}{i!} s(n,i, \mu, l)  &=&  \sum_{i=0}^{m-1} \dfrac{(-1)^i}{i!} s(m,i, \mu_0, l_0) \\
            &=& f_{m,1} = 0  \quad \quad (\text{because }m \text{ is even})
        \end{eqnarray*}
        \end{itemize}
        
        \item In this case, neither the leftmost nor the rightmost vertex is dominated by an element outside of the $m$-size subposet. Therefore : 
        \[s(n,i, \mu, l) = s(m,i, \mu_0, l_0)\]

        Hence : 
        \begin{eqnarray*}
            \sum_{i=0}^{m-1} \dfrac{(-1)^i}{i!} s(d,i, \mu, l)  &=&  \sum_{i=0}^{m-1} \dfrac{(-1)^i}{i!} s(m,i, \mu_0, l_0) \\
            &=& f_{1,m} = w_{-}(m) \quad \quad (\text{because }m \text{ is odd})
        \end{eqnarray*}

        \item In this case,$(\mu, l)$ can't be the result of a merge. Indeed, the right part and the left part are dominated by a vertex outside of the size-$m$ subposet. Therefore, neither of them can have the maximal label. Hence : 
        \[\sum_{i=0}^{m-1} \dfrac{(-1)^i}{i!} s(n,i, \mu, l)   = s(n,0,\mu, l) = \frac{((\frac{m-1}{2})!)^2}{m!} = w_{+}(m)\]
    \end{enumerate}

    \item we use the same argument as for the case $1$, when we are looking when is done the last mere. Hence :
    \begin{equation*}
        s(n,i, \mu_, l) =
            \begin{cases}
                  2 \displaystyle \sum_{\substack{j=1 \\ j \text{ odd} }}^{m-1} c_1(j) \times s(m-j, i-1,\mu_0, l_0)& \text{if } i \ge 1\\
                  c_1(m) & \text{if } i=0
            \end{cases}
        \end{equation*}

Therefore : 
\begin{eqnarray*}
             \sum_{i=0}^{m-1} \dfrac{(-1)^i}{i!} s(d,i,\mu, l)  &=&  c_1(m) + \sum_{i=1}^{m-1} 2\sum_{\substack{j=1 \\ j \text{ odd}}}^{m-1-i} \frac{(-1)^i}{i!} c_1(j) \times s(m-j, i-1,\mu_0, l_0) \\
             &=& c_1(m) + \sum_{\substack{j=1 \\ j \text{ odd}}}^{m-1} c_1(j)\sum_{i=0}^{m-1-j}\frac{(-1)^i}{i!}2s(m-j, i-1,\mu_0, l_0) \\
             &=& c_1(m) + \sum_{\substack{j=1 \\ j \text{ odd}}}^{m-1} c_1(j)f_{m-j,0} \\
             &=&  \sum_{\substack{j=1 \\ j \text{ odd}}}^{m} c_1(j)f_{m-j,0} \\   
\end{eqnarray*}\
The OGF of $f_{n,0}$ is $\sqrt{1-{z^2}/{4}}$   \\
The OGF of $c_1(j)$ with $j$ odd is : $ \dfrac{4\arcsin(\frac{z}{2})}{{2\sqrt{1-z^2/4}}}$. 
Hence, the result is the $m-$th coefficient of is $2\arcsin(z/2) = U(z)$
\end{proof}

\end{document}